\documentclass[12pt,reqno]{amsart}
\usepackage{latexsym}
\usepackage{amsmath}
\usepackage{amsfonts}
\usepackage{amssymb}

\usepackage{textcomp}

\usepackage{hyperref}

\def\={\equiv}
\renewcommand{\t}[1]{\widetilde{#1}}

\def\Z{\mathbb Z}

\vspace{0.5cm}

\author{Elena Kireeva, Alexei Krasilnikov, Vladimir Shchigolev}

\title[The finite basis problem for varieties modulo $x^3\equiv0$]{The finite basis property of varieties of associative algebras modulo $x^3\equiv0$}

\address[E.\;Kireeva]{Bauman Moscow State Technical University, ul. Baumanskaya 2-ya, 5/1, Moscow, 105005, Russia}
\email{eakireeva@yandex.ru}

\address[A.\;Krasilnikov]{Departamento de Matem\'atica, Universidade de Bras\'ilia, 70910-900 Bras\'ilia, DF, Brazil}
\email{alexei@unb.br}

\address[V.\;Shchigolev]{Financial University under the Government of the Russian Federation, Leningradsky Prospekt 49,
Moscow, 125993, Russia}
\email{shchigolev\_vladimir@yahoo.com}

\newtheorem{lem}{Lemma}
\newtheorem{teo}{Theorem}

\begin{document}

${}$

\vspace{-20pt}

\maketitle

\begin{abstract}
We prove that any variety of associative algebras over a Noetherian commutative ring
with unity satisfying the identity $x^3\=0$ is finitely based. This property has been known
for the varieties satisfying the identity $x^2\=0$. On the other hand, there are varieties
satisfying the identity $x^6\=0$ that are not finitely based.

Keywords: Associative algebras, polynomial identities, nil-alge\-bra, finite basis property.
\end{abstract}

\bigskip

\section{Introduction}

Let $K$ be a commutative associative ring with unity 
and $A$ be an associative algebra without unity over $K$ freely generated by (countably many) generators $x_1,x_2, \ldots$
We will often denote these generators by $x$ and $y$.
An element $v=v(x_1,\ldots,x_n)\in A$
is called a (polynomial) {\em  identity\/} of an associative $K$-algebra $M$ if
$v(f_1, \ldots , f_n) = 0$ for any $f_1, \ldots, f_n \in M$. In that case, we also say that $M$ {\it satisfies}
the identity $v\=0$. More generally, we say that $M$ {\it satisfies} an identity $v\=u$ if
$M$ satisfies the identity $v-u\=0$.

For example, any commutative algebra satisfies the identity $[x,y]\=0$, where $[x,y]=xy-yx$ is the commutator,
a nil-algebra satisfies the identity $x^n\=0$ for some positive integer $n$,
a nilpotent algebra of degree at most $s$ satisfies the identity $x_1 x_2 \cdots x_s\=0$.
If $\{v_t \mid t \in \Omega \}$ is a
set of identities, then the class of all $K$-algebras satisfying all
identities $v_t\=0$ for $t\in\Omega$ is called a {\em variety}.

A variety is {\em finitely based} if it can be defined by a finite set of identities.
It was proved by A.\,Kemer~\cite{kem} that all varieties of associative algebras over
fields of characteristic 0 are finitely based. On the other hand, A.\,Belov~\cite{belov}, A.\,Grishin~\cite{grish},
and V.\,Shchigolev~\cite{shch-1} constructed examples of non-finitely generated varieties
of associative algebras over fields of positive characteristic.
All these examples remain non-finitely based modulo $x^n\=0$ for $n$ big enough.

It is easy to see that any variety of associative algebras modulo $x^2\=0$ over a Noetherian ring
satisfies the identity $xy\=-yx$ and thus is finitely-based~\cite{kir}.
Moreover, according to the Dubnov-Ivanov-Nagata-Higman theorem~\cite{DI},~\cite{nag},~\cite{hig},
any variety of associative algebras modulo $x^3\=0$ over a field of characteristic $p$
with $p>3$ satisfies an identity of nilpotency. It is well known that modulo such an identity 
any variety of associative algebras is finitely based~\cite{bah}.
On the other hand, C.\,Gupta and A.\,Krasilnikov \cite{gup-kras} for algebras over fields of characteristic $2$
and E.\,Aladova and A.\,Krasilnikov \cite{al-kras} for algebras over fields of characteristic $3$
constructed examples of non-finitely based varieties of
associative algebras modulo $x^6\=0$.

Our main result is as follows.

\begin{teo} \label{teo-T-id}
Let $K$ be a Noetherian commutative associative ring with unity. Any variety of
associative $K$-algebras satisfying the identity
\begin{equation} \label{tozh-1}
x^3\=0
\end{equation}
is finitely based.
\end{teo}

An ideal $V$ of $A$ is called a {\em T-ideal,\/}
if $V$ is fully characteristic in $A$, that is, $\alpha (V) \subset V$ for any endomorphism $\alpha$ of
$A$. It is well known that there exists a one-to-one correspondence between
the set of all varieties of associative $K$-algebras and the set of all T-ideals of $A$.
This correspondence takes a variety ${\bf V}$ of associative algebras to the
T-ideal of all polynomial identities satisfied in all algebras of ${\bf V}$.
A variety of associative $K$-algebras is finitely based if and only if the corresponding T-ideal is finitely generated
as a T-ideal.

If ${\bf V}$ is a variety of associative $K$-algebras and $V$ is the corresponding T-ideal of $A$, then
the quotient algebra $A/V$ is called the {\em free algebra\/} of ${\bf V}$ or
a {\em relatively free algebra} freely generated by $x_1+V,x_2+V,\ldots$
A fully characteristic ideal in a relatively free algebra is also called a {\em T-ideal\/}.
There is a similar one-to-one correspondence between
the set of all subvarieties of ${\bf V}$ and the set of all T-ideals in $A/V$.

A submodule $U$ in a relatively free algebra $A/V$ is called a {\em T-space\/} if $U$ is fully characteristic
in $A/V$.
It was proved by V.\,Shchigolev~\cite{shch-3} for algebras over fields of characteristic 0
and by E.\,Kireeva~\cite{kir} for commutative algebras over Noetherian commutative associative rings with unity
that all T-spaces in the corresponding free algebras are finitely generated.
Also it is easy to see that all T-spaces in a nilpotent relatively free algebra are finitely generated.
On the other hand, V.Grishin \cite{grish} for algebras over fields of characteristic 2 and V.\,Shchigolev \cite{shch-2}
for algebras over fields of characteristic $p>2$ constructed examples of
non-finitely generated T-spaces.
All these examples remain non-finitely generated modulo $x^n\=0$ for some $n$.
In particular, in \cite{shch-2} examples of non-finitely generated T-spaces
modulo $x^3\=0$ in algebras over fields of characteristic 3 were constructed.
These T-spaces played an important role in
the negative solution of the finite basis problem in positive characteristic.

We can note that if all T-spaces in a relatively free algebra $A/V$ are finitely generated,
then all T-ideals in $A/V$ are finitely generated as well, and hence
all subvarieties of the variety $\bf{V}$ are finitely based.

Let $\mathbf{V_3}$ be the variety of associative $K$-algebras satisfying identity~(\ref{tozh-1})
and $V_3$ be the corresponding T-ideal in $A$. We will prove the following theorem.

\begin{teo} \label{teo-T-pr}
Let $K$ be a Noetherian commutative associative ring with unity
of characteristic 
not divisible by $3$.
Then all T-spaces in $A/V_3$ are finitely generated.
\end{teo}

\section{Auxiliary results}

This section serves as a preparation for the proofs of the theorems stated in the introduction.
It contains various calculations, some of them quite lengthy. So we encourage the reader to
use system of computer algebra (for example, Maxima) to verify them.
We abbreviate $\t{A}=A/V_3$.

\begin{lem}\label{lemma:1}
Let $K$ be a commutative associative ring with unity.
Then the algebra $\t{A}$ satisfies the identity
\begin{equation} \label{iden-2}
xyz+xzy+yxz+yzx+zxy+zyx\equiv0,
\end{equation}
which is the full linearization of~(\ref{tozh-1}).
\end{lem}
\begin{proof} By definition, $\t{A}$ satisfies the identity
\begin{multline*}
0\equiv(x+y+z)^3-(x+y)^3-(y+z)^3-(x+z)^3+x^3+y^3+z^3\\[6pt]
=xyz+xzy+yxz+yzx+zxy+zyx.
\end{multline*}
\end{proof}
\noindent
In what follows, we will denote the left-hand side of~(\ref{iden-2}) by $S_3(x,y,z)$.

\begin{lem}\label{lemma:3}
Let $K$ be an arbitrary commutative associative ring with unity of
characte\-ristic 2. Then $\t{A}$ satisfies the identity
\begin{equation} \label{tozh_p=2}
[x_1,x_2]x_3[x_4,x_5]\=0.
\end{equation}
\end{lem}
\begin{proof} Applying Lemma~\ref{lemma:1} and a routine (computer assisted) calculation, we get
that $\t{A}$ satisfies the identity
\begin{multline*}
0\=
S_3(x_3x_5x_1, x_4, x_2)+
	S_3(x_4, x_3, x_5x_1x_2)+
	S_3(x_3x_5x_4, x_1, x_2)\\
+
	S_3(x_2, x_1x_3x_4, x_5)+
	S_3(x_4, x_3, x_5x_2x_1)+
	S_3(x_2x_4, x_1, x_3x_5)\\
+
	S_3(x_4 x_1, x_5x_2,x_3)+
	S_3(x_1x_4x_5, x_2, x_3)+
	S_3(x_1x_5x_3, x_2, x_4)\\
+
	S_3(x_5, x_2x_3, x_4x_1)+
	S_3(x_5x_3x_1, x_4, x_2)+
	S_3(x_1x_2, x_4, x_5x_3)\\
+
	S_3(x_5, x_1x_2x_4, x_3)+
	S_3(x_5, x_4x_2, x_1x_3)+
	S_3(x_5, x_1x_4x_2, x_3)\\
+
	S_3(x_4x_5, x_2, x_3x_1)+
	S_3(x_4x_5x_2, x_3, x_1)+
	S_3(x_5, x_2x_1x_4, x_3)\\
+
	S_3(x_5, x_2x_4x_1, x_3)+
	S_3(x_2, x_1, x_4x_5x_3)+
	S_3(x_2, x_3x_1x_4, x_5)\\
+
	S_3(x_5x_3, x_1, x_4x_2)+
	S_3(x_2x_4, x_1, x_5x_3)+
	S_3(x_4x_1x_5, x_2, x_3)\\
+
	S_3(x_5, x_4x_1x_2, x_3)+
	S_3(x_2, x_3x_4, x_1x_5)+
	S_3(x_1, x_4, x_3x_5x_2)\\
+
	S_3(x_2, x_3x_4x_1, x_5)+
	S_3(x_5, x_4x_2x_1, x_3)+
	S_3(x_2x_3, x_4x_5, x_1)\\
+
	S_3(x_1, x_4, x_5x_2x_3)+
	S_3(x_2, x_1, x_5x_3x_4)+
	S_3(x_5x_3, x_4, x_2x_1)\\
+
	S_3(x_3x_5, x_2, x_1x_4)+
	S_3(x_3, x_5x_1, x_4x_2)+
	S_3(x_5x_1, x_3x_4, x_2)\\
=[x_1,x_2]x_3[x_4,x_5].
\end{multline*}
\end{proof}

\noindent
We will use left normed commutators inductively defined by the rules: $[x_1,x_2]=x_1x_2-x_2x_1$ (as in the introduction);
$[x_1,x_2,\ldots,x_n]=[[x_1,x_2,\ldots,x_{n-1}],x_n]$ for $n>2$.

\begin{lem}\label{lemma:4}
Let $K$ be an arbitrary commutative associative ring with unity of
characte\-ristic 3. Then $\t{A}$ satisfies the identity
\begin{equation} \label{tozh_p=3}
[x_1,x_2,x_3,x_4]=0.
\end{equation}
\end{lem}
\begin{proof}
Denoting $U(x_1,x_2,x_3)=[x_2,x_1,x_3]+[x_3,x_1,x_2]$, we get by Lemma~\ref{lemma:1} that $\t{A}$ satisfies
the identity
$$
U(x_1,x_2,x_3)\=U(x_1,x_2,x_3)+S_3(x_1,x_2,x_3)=3x_3x_1x_2+3x_2x_1x_3=0.
$$
Hence this algebra satisfies the identity
\begin{multline*}
0\=
U(x_1,x_2,x_3x_4)=x_3U(x_1,x_2,x_4)+U(x_1,x_2,x_3)x_4\\[4pt]
+[x_3,x_2][x_4,x_1]+[x_3,x_1][x_4,x_2]\=[x_3,x_2][x_4,x_1]+[x_3,x_1][x_4,x_2]
\end{multline*}
and thus also the identity
\begin{multline*}
[[x_1,x_2],[x_3,x_4]]=[x_1,x_2][x_3,x_4]-[x_3,x_4][x_1,x_2]\\[4pt]
\=[x_1,x_2][x_3,x_4]+[x_3,x_2][x_1,x_4]=[x_1,x_2][x_3,x_4]+[x_2,x_3][x_4,x_1]\\[4pt]
\=[x_1,x_2][x_3,x_4]-[x_2,x_1][x_4,x_3]=0.
\end{multline*}
Finally, as $U(x_3,[x_1,x_2],x_4)\=0$, we get
$$
[x_1,x_2,x_3,x_4]=[[x_1,x_2],x_3,x_4]
\=-[x_4,x_3,[x_1,x_2]]=[[x_4,x_3],[x_1,x_2]]\=0.
$$
\end{proof}

\begin{lem}\label{lemma:2}
Let $K$ be an arbitrary commutative associative ring with unity.
Then $\t{A}$ satisfies the identity
\begin{equation} \label{tozh_n-ne2,ne3}
2^4\,3^2\,x_1 x_2\cdots x_7\=0.
\end{equation}
\end{lem}
\begin{proof}
For $k=0,1,2,3$ we have
$$
0\=(x+ky)^3=x^3+k(x^2y+xyx+yx^2)+k^2(xy^2+yxy+y^2x)+k^3y^3
$$
``Solving'' this system by Cramer's rule, we get
$$
\left |
\begin{array}{cccc}
1&0&0&0\\
1&1&1&1\\
1&2&4&8\\
1&3&9&27
\end{array}
\right | (x^2y+xyx+yx^2)= 2^2\,3\,(x^2y+xyx+yx^2)\=0,
$$
We denote $g(x,y)=x^2y+xyx+yx^2$. Following \cite{lee}, we note that
$$
0\=2^2\,3\big( g(x,z)y^2+ g(x,zy)y+g(x,zy^2)
-xzg(y,x)-zg(y,x^2)\big)=2^2\,3^2\,x^2zy^2.
$$
Hence
\begin{multline*}
0\=2^2\,3^2\,x^2z(y+t)^2=2^2\,3^2\,x^2z(y^2+yt+ty+t^2)\\
=2^2\,3^2\,x^2zy^2+2^2 \, 3^2 \, x^2z(yt+ty)+2^2\,3^2\,x^2zt^2
\=2^2\,3^2\,x^2z(yt+ty).
\end{multline*}
Denoting $h(y,t)=yt+ty$, we get
$$
0\=2^2 \, 3^2 \, x^2z \left( h(y,u)t+h(y,ut)-uh(y,t) \right)=2^3 \, 3^2 \,x^2zyut,
$$
whence
\begin{multline*}
0\=2^3 \, 3^2 \, (x+v)^2zyut= 2^3 \, 3^2 \, (x^2+xv+vx+v^2)zyut=\\
=2^3 \, 3^2 \, x^2zyut + 2^3 \, 3^2 \, (xv+vx)zyut +2^3 \, 3^2 \, v^2zyut \\
\= 2^3 \, 3^2 \, h(x,v) \, zyut .
\end{multline*}
We conclude by
$$
0\=2^3\,3^2\,(h(x,w)v+h(x,wv)-wh(x,v))zyut = 2^4\,3^2\,xwvzyut.
$$
\end{proof}

\medskip

We also will use the following easy lemma, whose proof we leave to the reader.
\begin{lem} \label{mod-ch}
Let $I$, $U$, $V$ be submodules of a module $M$ such that
$U \subseteq V,\, U \cap \, I=V \cap\, I$ and $(U+I)/\,I=(V+I)/\,I$.
Then $U=V$.
\end{lem}

\section{Main proofs}

\begin{proof}[Proof of theorem \ref{teo-T-pr}]
{\it Part 1}. Let the characteristic of $K$ be equal to $p$, where $p$ is a prime number distinct from $3$.
If $p=2$ then by Lemma~\ref{lemma:3} the algebra $\t{A}$ satisfies identity (\ref{tozh_p=2}).
Substituting $[x_6,x_7]$ for $x_2$, we see that $\t{A}$ also satisfies the identity $[x_1,x_2][x_6,x_7][x_4,x_5]\=0$.
Thus by~\cite[Theorem 4]{kir-kras}, all T-spaces in $\t{A}$ are finitely generated.
If $p>3$ then $\t{A}$ satisfies identity (\ref{tozh_n-ne2,ne3}).
Since $2\cdot 1_K$ and $3\cdot 1_K$ are invertible in this case, the algebra $\t{A}$ is nilpotent and all T-spaces in $\t{A}$ are
finitely generated as well.

{\it Part 2}. Now let the characteristic of $K$ be equal to $p^s$, where $p$ is a prime number distinct from $3$ and
$s$ is a positive integer greater than $1$.
Let
\begin{equation} \label{T-sp-n-ch1}
U_1 \subseteq U_2 \subseteq \cdots
\end{equation}
be an increasing chain of T-spaces in $\t{A}$.
We are going to prove by decreasing induction on $r=0,\ldots,s$ that
this chain stabilizes if all its elements are contained in $p^r\t{A}$.
The base case $r=s$ is obvious, as $p^s\t{A}=0$.

Suppose that $r<s$ and $U_i\subseteq p^r\t{A}$ for any $i$.
We consider two chains
\begin{equation}\label{T-sp-n-ch2}
(U_1+p^{r+1}\t{A})/p^{r+1}\t{A} \subseteq (U_2+p^{r+1}\t{A})/p^{r+1}\t{A} \subseteq \cdots
\end{equation}
and
$$
U_1 \cap p^{r+1}\t{A} \subseteq U_2 \cap p^{r+1}\t{A} \subseteq \cdots
$$
By Lemma~\ref{mod-ch}, it suffices to prove that both chains stabilize.

The latter chain stabilizes by the inductive hypothesis.
To prove that the former chain also stabilizes, we consider the homomorphism
$$
\chi:\t{A}/p\t{A}\to p^r\t{A}/p^{r+1}\t{A}
$$
given by $f+p\t{A}\mapsto p^rf+p^{r+1}\t{A}$. It is clearly surjective.
As the modules $(U_i+p^r\t{A})/p^{r+1}\t{A}$ are T-spaces in $\t{A}/p^{r+1}\t{A}$,
their inverse images under $\chi$ are T-spaces in $\t{A}/p\t{A}$,
which is a relatively free algebra over $\Z/p\Z$. These inverse images form
an increasing chain, which stabilizes by Part~1. Applying $\chi$,
we get that chain~(\ref{T-sp-n-ch2}) also stabilizes.

{\it Part 3.} Now let the characteristic of $K$ be equal to a positive integer $n$ not divisible by 3.
We consider the factorization
$$
n=p_1^{s_1}p_2^{s_2}\cdots p_r^{s_r},
$$
where $p_1, \ldots, p_r$ are pairwise distinct prime numbers not equal to $3$.
We denote $n_i=n/p_i^{s_i}$ for all $i$. Since $n_1, \ldots, n_r$ are relatively prime,
there exist integers $u_1, \ldots, u_r$ such that $n_1u_1+ \cdots +n_ru_r=1$.
We denote $e_i=n_iu_i\cdot 1_K$. It is easy to check that $e_ie_j=0$ for $i \ne j$ and $e_i^2=e_i$.
Hence we get the following decomposition:
$$
\t{A} = e_1\t{A}\oplus e_2\t{A}\oplus\cdots\oplus e_r\t{A},
$$
where
$$
e_i\t{A}\cong\t{A}/p_i^{s_i}\t{A},
$$
The algebra $\t{A}/p_i^{s_i}\t{A}$ is an algebra over the ring $\Z/p_i^{s_i}\Z$.

Let
\begin{equation} \label{T-sp-n-ch8}
U_1 \subseteq U_2 \subseteq\cdots
\end{equation}
be an increasing chain of T-spaces in $\t{A}$.
Then we get the chain
$$
e_iU_1 \subseteq e_iU_2 \subseteq \cdots
$$
for each $i=1,\ldots,r$. Every such chain is a chain of T-spaces
in $e_i\t{A}$ and thus stabilizes by Part~2. Hence it follows that chain (\ref{T-sp-n-ch8})
stabilizes as well.
\end{proof}

\begin{proof}[Proof of Theorem \ref{teo-T-id}]
If $K$ is a ring of characteristic $3$, then by Lemma~\ref{lemma:4} the algebra $\t{A}$
satisfies identity (\ref{tozh_p=3}) and hence all T-ideals in $\t{A}$
are finitely generated by~\cite{popov}. Therefore, arguing similarly to the proof of Theorem~\ref{teo-T-pr},
we get that all T-ideals in $\t{A}$ are finitely generated if $K$ is a ring of characteristic $n\ne0$.

Let now $K$ be a ring of characteristic $0$. We denote by $I$ the T-ideal in $\t{A}$ generated by $x_1\cdots x_7+V_3$.
Let
\begin{equation} \label{T-id-0-ch1}
I_1 \subseteq I_2 \subseteq \cdots
\end{equation}
be an increasing chain of T-ideals in $\t{A}$.
Then we get two chains
$$
(I_1+I)/I \subseteq (I_2+I)/I \subseteq \cdots
$$
$$
I_1 \cap I \subseteq I_2 \cap I \subseteq \cdots
$$
By Lemma~\ref{mod-ch}, it suffices to prove that both chains stabilize.

The ideals $(I_j +I)/I$ are T-ideals in the nilpotent algebra $\t{A}/I$ and hence the former chain 
stabilizes.

To prove that the latter chain stabilizes, we consider the T-ideal $N=2^4 3^2\t{A}$ in $\t{A}$.
Again it suffices to prove that both chains
\begin{equation} \label{T-id-0-ch4}
((I_1 \cap I)+ N)/N \subseteq ((I_2 \cap I)+N)/N \subseteq\cdots
\end{equation}
\begin{equation} \label{T-id-0-ch5}
I_1 \cap I \cap N \subseteq I_2 \cap I \cap N \subseteq\cdots
\end{equation}
stabilize.

The ideals $((I_j \cap I)+ N)/N$ are T-ideals in the algebra $\t{A}/N$, which is an algebra
over $\Z/2^4 3^2\Z$ and hence chain (\ref{T-id-0-ch4}) stabilizes.
Chain (\ref{T-id-0-ch5}) stabilizes since $I \cap N$ is the T-ideal in $\t{A}$ generated by
$2^4 3^2 x_1 \ldots x_7+V_3$ and hence $I \cap N=\{0\}$ by Lemma~\ref{lemma:2}.
\end{proof}

\def\sp{\vspace{4pt}}

\end{document}